\documentclass{amsart}

\usepackage{fancyhdr}
\usepackage{lastpage}
\usepackage{stmaryrd,yhmath}

\newcommand{\bdis}{\begin{displaymath}}
\newcommand{\edis}{\end{displaymath}}
\newcommand{\be}{\begin{equation}}
\newcommand{\ee}{\end{equation}}
\newcommand{\mbb}{\mathbb}
\newcommand{\mcal}{\mathcal}

\newcommand{\vp}{\varphi}

\newcommand{\vth}{\vartheta}

\newcommand{\zf}{\zeta\left(\frac{1}{2}+it\right)}

\theoremstyle{definition}

\theoremstyle{remark}
\newtheorem{remark}[]{Remark}

\newtheorem*{mydef1}{{\bf Theorem}}

\newtheorem*{mydef4}{{\bf Corollary}}

\newtheorem*{mydef7}{{\bf Question}}

\numberwithin{equation}{section}

\begin{document}

\title{Jacob's ladders, $\zeta$-factorization and infinite set of metamorphosis of a multiform}

\author{Jan Moser}

\address{Department of Mathematical Analysis and Numerical Mathematics, Comenius University, Mlynska Dolina M105, 842 48 Bratislava, SLOVAKIA}

\email{jan.mozer@fmph.uniba.sk}


\begin{abstract}
In this paper we use Jacob's ladders together with fundamental Hardy-Littlewood formula (1921) to prove
the so-called $\zeta$-factorization formula on the critical line. Simultaneously, we obtain a set of control parameters
of metamorphosis of a multiform connected with the Riemann-Siegel formula.
\end{abstract}
\maketitle

\section{Introduction and main result}

\subsection{}

Let us remind the unique factorization theorem: any positive integer ($\not= 1$) can be expressed as a product of primes. This expression is unique except
of the order in which the primes occur. Thus, we have
\be \label{1.1}
\begin{split}
& n=p_1^{m_1}p_2^{m_2}\cdots p_r^{m_r},\ n\in\mbb{N},\ n\not=1, \\
& p_1<p_2<\dots<p_r,\ m_1,m_2,\dots,m_r\in\mbb{N},
\end{split}
\ee
where
\bdis
p_1,\dots,p_r
\edis
are primes. Next, Euler's identity (1737)
\be \label{1.2}
\sum_{n=1}^\infty \frac{1}{n^x}=\prod_p (1-p^{-x})^{-1},\ x>1
\ee
may be regarded as an analytical equivalent of the unique factorization theorem.

Further, Riemann has defined in 1859 the zeta-function by the formula
\be \label{1.3}
\zeta(s)=\sum_{n=1}^\infty \frac{1}{n^s}=\prod_p (1-p^{-s})^{-1},\ s=\sigma+it\in\mbb{C},\ \sigma>1,
\ee
and, consequently, by the method of analytic continuation, he defined the function $\zeta(s)$ at all
finite points $s\in\mbb{C}$ except for a simple pole at $s=1$. Now, since from (\ref{1.1})), (comp.
(\ref{1.2}), (\ref{1.3})) the exponential form of the factorization
\be \label{1.4}
\left(\frac{1}{n}\right)^x=\prod_{l=1}^r\left(\frac{1}{p_l^{m_l}}\right)^x,\ x>1,
\ee
follows (as an example), then we can formulate the following.

\begin{remark}
On the basis of the analytic continuation of the Riemann's zeta-function we may suppose that:
\begin{itemize}
\item[(a)] there are inherited forms of the factorization (\ref{1.4}) in the region $\sigma\leq 1,\ s\not=1$, (together
with mutation of these),
\item[(b)] especially, we may suppose that on the critical line itself
$$\sigma=1/2$$
there are some new species of the factorization property (\ref{1.4}).
\end{itemize}
\end{remark}

\subsection{}

In this direction  the following theorem holds true.

\begin{mydef1}
Let
\bdis
k=1,\dots,k_0\in \mbb{N}
\edis
for every fixed $k_0$ and
\be \label{1.5}
H=H(T)\in \left(\frac{\ln\ln T}{\ln T},T^{\frac{1}{\ln\ln T}}\right)
\ee
for every sufficiently big $T>0$. Then for every $k$ and for every $H$ that fulfils condition (\ref{1.5}) there are the functions
\be \label{1.6}
\begin{split}
& H_k=H_k(T,H)>0 , \\
& \alpha_r=\alpha_r(T,H,k)>0,\ r=0,1,\dots,k , \\
& \alpha_r\not=\gamma:\ \zeta\left(\frac 12+i\gamma\right)=0,
\end{split}
\ee
such that the following $\zeta$-factorization formula
\be \label{1.7}
\begin{split}
& \sqrt{\frac{\Lambda}{\left|\zeta\left(\frac 12+i\alpha_0\right)\right|}}\sim \prod_{r=1}^k \left|\zeta\left(\frac 12+i\alpha_r\right)\right|, \\
& \Lambda=\Lambda(T,H,k)=\sqrt{2\pi}\frac{\sqrt{H}}{H_k}\ln^kT
\end{split}
\ee
holds true. Moreover, the sequence
\bdis
\{ \alpha_r\}_{r=0}^k
\edis
has the following properties:
\be \label{1.8}
T<\alpha_0<\alpha_1<\dots<\alpha_k ,
\ee
\be \label{1.9}
\alpha_{r+1}-\alpha_r\sim (1-c)\pi(T),\ r=0,1,\dots,k-1,
\ee
where
\bdis
\pi(T)\sim \frac{T}{\ln T}, \ T\to\infty
\edis
is the prime counting function and $c$ is the Euler's constant.
\end{mydef1}

\begin{remark}
The asymptotic behavior of the set
\be \label{1.10}
\{ \alpha_0,\alpha_1,\dots\alpha_k\}
\ee
by (\ref{1.9}) is as follows: if $T\to\infty$ then the points of the set (\ref{1.10}) recede unboundedly each from other and all these points together
recede to infinity. Hence, at $T\to\infty$ the set (\ref{1.10}) behaves like one dimensional Friedmann-Hubble expanding universe.
\end{remark}

\begin{remark}
Consequently, the following holds true: even if the distances
\bdis
\begin{split}
& |\alpha_q-\alpha_r|\sim |q-r|(1-c)\pi(T)\to\infty \ \mbox{as}\ T\to\infty,  \\
& q,r=0,1,\dots,k,\ q\not=r
\end{split}
\edis
(see (\ref{1.9})) of elements of the set (\ref{1.10}) are gigantic, there is a close constraint between the values of the set
\bdis
\left\{\left|\zeta\left(\frac 12+i\alpha_0\right)\right|,\left|\zeta\left(\frac 12+i\alpha_1\right)\right|,\dots,\left|\zeta\left(\frac 12+i\alpha_k\right)\right|\right\},
\edis
namely the factorization formula (\ref{1.7}).
\end{remark}

\section{On metamorphosis of a multiform generated by the Riemann-Siegel formula}

\subsection{}

Let us remind the Riemann-Siegel formula
\be \label{2.1} \begin{split}
& Z(t)=2\sum_{n\leq\tau(t)}\frac{1}{\sqrt{n}}\cos\{\vth(t)-t\ln n\}+\mcal{O}(t^{-1/4}), \\
& \tau(t)=\sqrt{\frac{t}{2\pi}},
\end{split}
\ee
(see \cite{6}, p. 60, comp. \cite{7}, p. 79) where
\bdis
\vth(t)=-\frac t2\ln\pi+\mbox{Im}\ln\Gamma\left(\frac 14+i\frac 12\right),
\edis
(see (\cite{7}), p. 239). Next, we put in (\ref{2.1})
\be \label{2.2}
Z(t)=\sum_{n\leq\tau(t)}a_nf_n(t)+R(t),
\ee
where
\be \label{2.3} \begin{split}
& a_n= \frac{1}{\sqrt{n}}, \\
& f_n(t)=\cos\{\vth(t)-t\ln n\}, \\
& R(t)=\mcal{O}(t^{-1/4}),
\end{split}
\ee
and the functions $f_n(t)$ are nonlinear. That is, the function $Z(t)$ is nonlinear monoform. Now, we define the
following multiform
\be \label{2.4}
G(x_1,\dots,x_k)=\prod_{r=1}^k |Z(x_r)|,\quad x_r>T>0,\ k\geq 2,
\ee
where the monoform $Z(t)$ is the generating function of $G$.

Further, we define the subset of all random sample points
\bdis
\begin{split}
& (\bar{x}_1,\dots,\bar{x}_k):\ T<\bar{x}_1<\bar{x}_2<\dots<\bar{x}_k, \\
& \bar{x}_r\not=\gamma:\ \zeta\left(\frac 12+i\gamma\right)=0,\ r=1,\dots,k ,
\end{split}
\edis
and propose the following

\begin{mydef7}
Is there in this subset a point
\bdis
(\bar{y}_1,\dots,\bar{y}_k)
\edis
of metamorphosis of the multiform (\ref{2.4}) (that is, the point of significant change of the structure of
the multiform (\ref{2.4}))?
\end{mydef7}

The answer is given in the following.

\begin{mydef4}
\be \label{2.5}
\begin{split}
& \prod_{r=1}^k \left|\sum_{n\leq\tau(\alpha_r)}a_nf_n(\alpha_r)+R(\alpha_r)\right|\sim \\
& \sim \sqrt{\frac{\Lambda}{\left|\sum_{n\leq\tau(\alpha_0)}a_nf_n(\alpha_0)+R(\alpha_0)\right|}},\ T\to\infty,
\end{split}
\ee
i.e. to the infinite subset of the points
\bdis
\{\alpha_1(T),\alpha_2(T),\dots,\alpha_k(T)\},\quad T\in (T_0,\infty),
\edis
where $T_0$ is a sufficiently big, an infinite set of metamorphoses of the multiform (\ref{2.4}) into quite distinct form on the right-hand side of (\ref{2.5}) corresponds.
\end{mydef4}

\begin{remark}
We shall call the elements of the set
\be \label{2.6}
\{\alpha_0(T),\alpha_1(T),\dots,\alpha_k(T)\},\ T\in (T_0,+\infty)
\ee
as \emph{control parameters} (\emph{functions}) of metamorphosis. The reason to this is that the parameters
\bdis
\alpha_1(T),\dots,\alpha_k(T)
\edis
change the old form into the new one (see (\ref{2.5}), and this last is controlled by the parameter $\alpha_0(T)$.
That is, the set (\ref{2.6}) plays a similar role as shem-ha-m'forash in Golem's metamorphosis.
\end{remark}

\section{On Hardy-Littlewood fundamental Lemma 18 from the memoir \cite{2}}

\subsection{}

Let us remind the following (see \cite{2}, pp. 304, 305): if
\be \label{3.1}
\begin{split}
& J=J(T,H)=\int_T^{T+U}\mcal{J}^2{\rm d}t, \\
& \mcal{J}=\mcal{J}(t,H)=\int_t^{t+H} x(u){\rm d}u,
\end{split}
\ee
then the following Hardy-Littlewood formula
\be \label{3.2}
J=\int_T^{T+H}\mcal{J}^2{\rm d}t=\pi\sqrt{2\pi}HU+\mcal{O}\left(\frac{U}{\ln T}\right)
\ee
holds true, where
\be \label{3.3}
U\in [T^a,T^b],\ \frac 12<a<b\leq \frac 58,\ 0<H\leq T^{\epsilon},
\ee
and $\epsilon$ is positive and sufficiently small, and (see \cite{2}, p. 290)
\be \label{3.4}
\zf=-\left(\frac 2\pi\right)^{1/4}e^{i\pi/8}(2\pi e)^{\frac 12 it}e^{-\frac 12 it\ln t}x(t)\left\{ 1+\mcal{O}\left(\frac 1t\right)\right\}.
\ee

\subsection{}

We use the Riemann function $Z(t)$ instead of the Hardy-Littlewood function $x(t)$. Namely, the formula
\be \label{3.5}
\zf=e^{-i\vth(t)}Z(t)
\ee
together with formula for $\vth(t)$ from \cite{7}, p. 329 give that
\be \label{3.6}
x(t)=-\left(\frac{\pi}{2}\right)^{1/4}Z(t)\left\{ 1+\mcal{O}\left(\frac 1t\right)\right\}.
\ee
Hence, we use the following variant of the Hardy-Littlewood formula (see (\ref{3.6}))
\be \label{3.7}
\bar{J}=\bar{J}(T,H)=\int_T^{T+U_0}\bar{\mcal{J}}^2{\rm d}t=2\pi HU_0+\mcal{O}\left(\frac{U_0}{\ln T}\right),
\ee
where (comp. (\ref{3.1}) -- (\ref{3.3}))
\be \label{3.8}
\begin{split}
& \bar{\mcal{J}}=\bar{\mcal{J}}(t,H)=\int_t^{t+H}Z(u){\rm d}u, \\
& U_0=T^{0.5001}, \\
& H\in \left(\frac{\ln\ln T}{\ln T},T^{\frac{1}{\ln\ln T}}\right).
\end{split}
\ee

\begin{remark}
Since (see (\ref{3.7}))
\bdis
2\pi HU_0+\mcal{O}\left(\frac{U_0}{\ln T}\right)=2\pi HU_0\left\{  1+\mcal{O}\left(\frac{1}{H\ln T}\right)\right\},
\edis
then we have that the formula (\ref{3.7}) is asymptotic formula for $H$ of (\ref{3.8}).
\end{remark}

\subsection{}

Let us remind the following sentences of Hardy and Littlewood (see \cite{2}, p. 315): \emph{As was observed in 5.5, we do not use
the full force of Lemma 18. The complete lemma, however, seems of considerable interest in itself, and it may prove to be of service
in the future}.

\begin{remark}
We notice explicitly, that our theory of the Jacob's ladders together with the Hardy-Littlewood asymptotic formula (\ref{3.7}) constitute
the basis of our result about metamorphosis of corresponding multiform. Thus, we have obtained, after 94 years , the result that use the full force
of the asymptotic Hardy-Littlewood formula from their Lemma 18.
\end{remark}

\section{Jacob's ladders and the Hardy-Littlewood integral (1918)}

Let us remind that we have introduced (see \cite{4}, (9.1), (9.2)) the following formula
\be \label{4.1}
\tilde{Z}^2(t)=\frac{{\rm d}\vp_1(t)}{{\rm d}t},
\ee
where
\be \label{4.2}
\begin{split}
& \tilde{Z}^2(t)=\frac{Z^2(t)}{2\Phi'_{\vp}[\vp(t)]}=\frac{\left|\zf\right|^2}{\omega(t)}, \\
& \omega(t)=\left\{ 1+\mcal{O}\left(\frac{\ln\ln t}{\ln t}\right)\right\}\ln t.
\end{split}
\ee
The function
\bdis
\vp_1(t)
\edis
is called the Jacob's ladder (see our paper \cite{3}) according to Jacob's dream in Chumash, Bereishis, 28:12, has the following properties:
\begin{itemize}
\item[(a)]
\bdis
\vp_1(t)=\frac 12\vp(t),
\edis
\item[(b)] the function $\vp(t)$ is a solution of the nonlinear integral equation (see \cite{3}, \cite{4})
\bdis
\int_0^{\mu[x(T)]}Z^2(t)e^{-\frac{2}{x(T)}t}{\rm d}t=\int_0^T Z^2(t){\rm d}t,
\edis
where each admissible function $\mu(y)$ generates the solution
\bdis
y=\vp_\mu(T)=\vp(T),\ \mu(T)\geq 7y\ln y.
\edis
\end{itemize}

\begin{remark}
The main goal of introducing the Jacob's ladders is described in \cite{3}, where we have shown, by making use of these Jacob's ladders, that the
Hardy-Littlewood integral (1918)
\bdis
\int_0^T \left|\zf\right|^2{\rm d}t
\edis
(see \cite{1}, pp. 122, 151 -- 156) has - in addition to the Hardy-Littlewood expression (and also other similar to this one) possessing
an unbounded error at $T\to\infty$ - the following infinite set of almost exact expressions
\bdis
\begin{split}
 & \int_0^T \left|\zf\right|^2{\rm d}t=\vp_1(T)\ln\vp_1(T)+(c-\ln 2\pi)\vp_1(T)+ \\
 & + c_0+\mcal{O}\left(\frac{\ln T}{T}\right),\quad T\to\infty,
\end{split}
\edis
where $c_0$ is the constant from the Titchmarsh-Kober-Atkinson formula (see \cite{7}, p. 114).
\end{remark}

\begin{remark}
 The Jacob's ladder $\vp_1(T)$ can be interpreted by our formula (see \cite{3}, (6.2))
 \bdis
 T-\vp_1(T)\sim (1-c)\pi(T),
 \edis
 where $\pi(T)$ is the prime-counting function, as an asymptotic complementary function to the
 function
 \bdis
 (1-c)\pi(T)
 \edis
 in the following sense
 \bdis
 \vp_1(T)+(1-c)\pi(T)\sim T,\quad T\to\infty.
 \edis
\end{remark}

\section{Proof of Theorem}

\subsection{}

First of all we obtain from the formula
\bdis
\int_T^{T+U_0}\left\{\int_t^{t+H}Z(u){\rm d}u\right\}^2{\rm d}t\sim 2\pi HU_0
\edis
(see (\ref{3.7}), (\ref{3.8})), by mean-value theorem, that
\bdis
\left\{\int_\eta^{\eta+H} Z(u){\rm d}u\right\}^2\sim 2\pi H,
\edis
i.e.
\be \label{5.1}
\left| \int_\eta^{\eta+H} Z(u){\rm d}u\right|\sim \sqrt{2\pi H},\ \eta=\eta(T,H)\in
(T,T+U_0).
\ee
Since (see (\ref{3.8}))
\be \label{5.2}
H=o\left(\frac{T}{\ln T}\right),
\ee
then by our lemma (comp. \cite{5}, (7.1), (7.2))
\bdis
\left|\int_{\overset{k}{\eta}}^{\overset{k}{\wideparen{\eta+H}}}Z[\vp_1^k(t)]
\prod_{r=0}^{k-1}\tilde{Z}^2[\vp_1^r(t)]{\rm d}t\right|\sim \sqrt{2\pi H}
\edis
and, by the mean-value theorem, we obtain that
\be \label{5.3}
\begin{split}
 & |Z[\vp_1^k(\beta)]|\sim \sqrt{2\pi}\frac{\sqrt{H}}{\overset{k}{\wideparen{\eta+H}}-\overset{k}{\eta}}
 \prod_{r=0}^{k-1}\tilde{Z}^{-2}[\vp_1^r(\beta)], \\
 & \beta=\beta(T,H,k)\in (\overset{k}{\eta},\overset{k}{\wideparen{\eta+H}}).
\end{split}
\ee

\subsection{}

Now, we make a transformation of the formula (\ref{5.3}). First of all, we have (comp. \cite{5},
Property 2, (6.4))
\bdis
\beta\in (\overset{k}{\eta},\overset{k}{\wideparen{\eta+H}}) \ \Rightarrow \
\vp_1^r(\beta)\in (\overset{k-r}{\eta},\overset{k-r}{\wideparen{\eta+H}}),\
r=0,1,\dots,k,
\edis
i.e.
\be \label{5.4}
\begin{split}
 & \vp_1^0(\beta)=\alpha_k\in (\overset{k}{\eta},\overset{k}{\wideparen{\eta+H}}), \\
 & \vp_1^1(\beta)=\alpha_{k-1}\in (\overset{k-1}{\eta},\overset{k-1}{\wideparen{\eta+H}}), \\
 & \vdots \\
 & \vp_1^{k-2}(\beta)=\alpha_2\in (\overset{2}{\eta},\overset{2}{\wideparen{\eta+H}}), \\
 & \vp_1^{k-1}(\beta)=\alpha_1\in (\overset{1}{\eta},\overset{1}{\wideparen{\eta+H}}), \\
 & \vp_1^k(\beta)=\alpha_0\in (\eta,\eta+H).
\end{split}
\ee
Consequently, from (\ref{5.3}) by (\ref{5.4}), (\ref{4.2}) the formula
\be \label{5.5}
\begin{split}
 & \left|\zeta\left(\frac 12+i\alpha_0\right)\right|\sim \sqrt{2\pi}
 \frac{\sqrt{H}}{H_k}\prod_{l=1}^k \frac{\omega(\alpha_l)}
 {\left|\zeta\left(\frac 12+i\alpha_l\right)\right|^2}, \\
 & \alpha_r=\alpha_r(T,H,k),\ r=0,1,\dots,k,\ H_k(T,H)=\overset{k}{\wideparen{\eta+H}}-\overset{k}{\eta}
\end{split}
\ee
follows.

\subsection{}

Next, let us remind the following properties of the disconnected set
\bdis
\Delta(T,H,k)=\bigcup_{r=0}^k [\overset{r}{\eta},\overset{r}{\wideparen{\eta+H}}].
\edis
(see (\ref{5.2}), comp. \cite{5}, (2.5) -- (2.7), (2.9)) . If
\bdis
H=o\left(\frac{T}{\ln T}\right),
\edis
then
\be \label{5.6}
\begin{split}
 & |[\overset{r}{\eta},\overset{r}{\wideparen{\eta+H}}]|=\overset{r}{\wideparen{\eta+H}}-\overset{r}{\eta}=
 o\left(\frac{T}{\ln T}\right),\ r=1,\dots,k , \\
 & |[\overset{r-1}{\wideparen{\eta+H}},\overset{r}{\eta}]|\sim (1-c)\pi(T), \\
 & [\eta,\eta+H]\prec [\overset{1}{\eta},\overset{1}{\wideparen{\eta+H}}]\prec \dots \prec
 [\overset{k}{\eta},\overset{k}{\wideparen{\eta+H}}].
\end{split}
\ee
Hence, from (\ref{5.4}), (\ref{5.6}) the properties (\ref{1.8}), (\ref{1.9}) follow immediately. Further,
we have (comp. \cite{5}, (4.3))
\bdis
\ln t\sim \ln \eta,\ \forall\- t\in (\eta,\overset{k}{\wideparen{\eta+H}}),
\edis
and (see (\ref{3.8}), (\ref{5.1}))
\be \label{5.7}
\ln t\sim \ln\eta\sim\ln T,\ \forall\- t\in (\eta,\overset{k}{\wideparen{\eta+H}}) .
\ee
Consequently, we obtain from (\ref{5.5}) by (\ref{4.2}), (\ref{5.7}) the following
\bdis
\left|\zeta\left(\frac 12+i\alpha_0\right)\right|\sim \sqrt{2\pi}\frac{\sqrt{H}}{H_k}
\ln^kT\prod_{r=1}^k \left|\zeta\left(\frac 12+i\alpha_r\right)\right|^{-2},
\edis
i.e. the formula (\ref{1.7}) holds true.

\section{Concluding remarks}

If we use the formulae (comp. \cite{7}, pp. 221, 329)
\bdis
\begin{split}
 & \vth(t)=\frac t2\ln\frac{t}{2\pi}-\frac t2-\frac{\pi}{8}+\mcal{O}\left(\frac 1t\right), \\
 & \vth'(t)=\frac 12\ln\frac{t}{2\pi}+\mcal{O}\left(\frac 1t\right), \\
 & \vth''(t)\sim \frac{1}{2t}, \\
 & t\to\infty
\end{split}
\edis
then we obtain from (\ref{2.1}), by a little of transformations, the local kind of the spectral
representation
\be \label{6.1}
\begin{split}
 & Z(x_r)=2\sum_{n\leq\tau(x_r)}\frac{1}{\sqrt{n}}\cos\left\{ t\ln\frac{\tau(x_r)}{n}-\frac{x_r}{2}-
 \frac{\pi}{8}\right\}+\\
 & + \mcal{O}(x_r^{-1/4}), \\
 & t\in [x_r,x_r+H], H\in \left(\frac{\ln\ln T}{\ln T},T^{\frac{A}{\ln\ln T}}\right),\
 \tau(x_r)=\sqrt{\frac{x_r}{2\pi}},
\end{split}
\ee
i.e. the Riemann-Siegel formula (\ref{2.1}), $t=x_r$.

\begin{remark}
 Namely, the sequence
 \be \label{6.2}
 \{\omega_{n,r}\}_{n\leq \tau(x_r)},\ \omega_{n,r}=\ln\frac{\tau(x_r)}{n},\ r=1,\dots,k
 \ee
 of the cyclic frequencies $\omega_{n,r}$ will be called as the local spectrum of the
 Riemann-Siegel formula (\ref{2.1}).
\end{remark}

\begin{remark}
Consequently, we have (see (\ref{2.4}), (\ref{6.1}), (\ref{6.2})) that the multiform
\bdis
G(x_1,x_2,\dots,x_k)
\edis
expresses the complicated oscillating process.Just for this oscillating multiform we have constructed the set
of metamorphosis described by the formula (\ref{2.5}).
\end{remark}

\thanks{I would like to thank Michal Demetrian for his help with electronic version of this paper.}

\end{document}